\documentclass[11pt]{amsart}
\usepackage[T1]{fontenc}
\usepackage{lmodern}
\usepackage{amsmath,amssymb}
\usepackage{mathtools}
\usepackage[margin=1.1in]{geometry}
\usepackage{tikz}
\usetikzlibrary{matrix,calc}
\usepackage{amsmath, amssymb, amsthm}
\usepackage{cite}
\usepackage{enumitem}
\usepackage{comment}

\usepackage{amsmath}
\usepackage{amsfonts}
\usepackage{tikz}
\usepackage{tikz-cd}

\theoremstyle{plain}
\newtheorem{proposition}{Proposition}
\newtheorem{thm}[proposition]{Theorem}
\newtheorem{lemma}[proposition]{Lemma}
\newtheorem{conjecture}[proposition]{Conjecture}

\theoremstyle{definition}
\newtheorem{definition}[proposition]{Definition}
\newtheorem{example}[proposition]{Example}

\theoremstyle{remark}
\newtheorem{remark}[proposition]{Remark}

\title{Coxeter Group Symmetries of Frieze Patterns}
\author{Yuta Kozaki}
\address{Department of Mathematics, Graduate School of Science, Rikkyo
University, Toshima-ku, Tokyo 171-8501, Japan.}
\email{26rc001e@rikkyo.ac.jp}
\author{Saburo Kakei}
\address{Department of Mathematics, College of Science, Rikkyo University, Toshima-ku, Tokyo 171-8501, Japan.} 
\email{kakei@rikkyo.ac.jp}
\begin{document}

\maketitle

\begin{abstract}
Coxeter group symmetries of frieze patterns 
are investigated through their relation to $2\times n$ matrices
and the associated Pl\"ucker coordinates. 
For odd order $n$, the symmetric group  
acts on the set of friezes, while for even $n$, 
the dihedral subgroup acts. 
We further show that 
a larger group acts when $n = 6$.
\end{abstract}

\section{Introduction}

The study of frieze patterns was initiated by Coxeter \cite{Coxeter1971}.
In his original formulation, a frieze pattern is represented as an infinite
planar array
of numbers arranged in staggered rows.
%A frieze of order $n$ (which corresponds to a width of $n-3$) is constructed according to the following two rules. The first
\begin{comment}
\begin{enumerate}[label=(\roman*)]
    \item \textbf{Boundary conditions:}
    The pattern consists of $n+1$ horizontal rows.
    To define the boundary of the frieze, the first and last rows (Row 0 and Row $n$)　are composed entirely of zeros, while the second and second-to-last rows (Row 1 and Row $n-1$) consist entirely of ones. In other words, for any integer $i$, we have:
    \begin{equation}
    f_{i,i} = f_{i,i+n} = 0, \quad f_{i,i+1} = f_{i,i+n-1} = 1.
    \end{equation}

    \item \textbf{The unimodular rule:}
    The internal entries of the frieze are governed by a local arithmetic property known as the ``unimodular rule''. Given any four adjacent elements arranged in a diamond configuration:
    \[
    \begin{array}{ccc}
          & f_{i+1,j} &   \\
        f_{i,j} &   & f_{i+1,j+1} \\
          & f_{i,j+1} &
    \end{array}
    \]
    the following relation must hold:
    \begin{equation}\label{eq:diamond-determinant}
        f_{i,j}f_{i+1,j+1} - f_{i+1,j}f_{i,j+1} = 1.
    \end{equation}
    for $i<j<i+n$.
\end{enumerate}
\end{comment}
A frieze pattern of order $n$ is a grid of numbers with $n+1$
rows as in Figure \ref{fig.1},
where the first and last rows (row $0$ and row $n$) are
composed entirely of zeros,
while the second and second-to-last rows (row $1$ and row $n-1$)
consist entirely of ones.
In other words, for any integer $i$, we have
    \begin{equation}
    f_{i,i} = f_{i,i+n} = 0, \quad f_{i,i+1} = f_{i,i+n-1} = 1.
    \end{equation}
The internal entries of the frieze are governed by the ``unimodular rule'';
given any four adjacent elements arranged in a diamond configuration,
\[
    \begin{array}{ccc}
          & f_{i+1,j} &   \\
        f_{i,j} &   & f_{i+1,j+1}, \\
          & f_{i,j+1} &
    \end{array}
    \]
    the relation
    \begin{equation}\label{eq:diamond-determinant}
        f_{i,j}f_{i+1,j+1} - f_{i+1,j}f_{i,j+1} = 1
    \end{equation}
    must hold for $i<j<i+n$.
By the unimodular rule \eqref{eq:diamond-determinant},
providing $n-3$ diagonal entries $f_{1,3}, f_{1,4}, \dots, f_{1,n-1}$ is sufficient to determine
all other entries \cite{Coxeter1971}.

\begin{figure}[htbp]
\centering
\begin{tikzpicture}
    \draw (-0.5,0.0) -- (-0.5,-4.0);

    \node at (-1.25,0.25) {row};
    \node at (-1.25,-0.25) {0};
    \node at (-1.25,-0.75) {1};
    \node at (-1.25,-1.25) {2};
    \node at (-1.25,-1.75) {3};
    \node at (-1.25,-2.75) {$n-2$};
    \node at (-1.25,-3.25) {$n-1$};
    \node at (-1.25,-3.75) {$n$};

    \node at (0.5,-0.25) {0};
    \node at (2.0,-0.25) {0};
    \node at (3.5,-0.25) {0};
    \node at (5.0,-0.25) {0};
    \node at (6.5,-0.25) {0};

    \node at (1.25,-0.75) {1};
    \node at (2.75,-0.75) {1};
    \node at (4.25,-0.75) {1};
    \node at (5.75,-0.75) {1};
    \node at (7.25,-0.75) {1};

    \node at (2.0,-1.25) {$f_{1,3}$};
    \node at (3.5,-1.25) {$f_{2,4}$};
    \node at (5.0,-1.25) {$f_{3,5}$};
    \node at (6.5,-1.25) {$f_{4,6}$};
    \node at (8.0,-1.25) {$f_{5,7}$};

    \node at (2.75,-1.75) {$f_{1,4}$};
    \node at (4.25,-1.75) {$f_{2,5}$};
    \node at (5.75,-1.75) {$f_{3,6}$};
    \node at (7.25,-1.75) {$f_{4,7}$};
    \node at (8.75,-1.75) {$f_{5,8}$};

    \node at (3.5,-2.25) {$\ddots$};
    \node at (5.0,-2.25) {$\ddots$};
    \node at (6.5,-2.25) {$\ddots$};
    \node at (8.0,-2.25) {$\ddots$};
    \node at (9.5,-2.25) {$\ddots$};

    \node at (4.25,-2.75) {$f_{1,n-1}$};
    \node at (5.75,-2.75) {$f_{2,n}$};
    \node at (7.25,-2.75) {$f_{3,n+1}$};
    \node at (8.75,-2.75) {$f_{4,n+2}$};
    \node at (10.25,-2.75) {$f_{5,n+3}$};

    \node at (5.0,-3.25) {1};
    \node at (6.5,-3.25) {1};
    \node at (8.0,-3.25) {1};
    \node at (9.5,-3.25) {1};
    \node at (11.0,-3.25) {1};

    \node at (5.75,-3.75) {0};
    \node at (7.25,-3.75) {0};
    \node at (8.75,-3.75) {0};
    \node at (10.25,-3.75) {0};
    \node at (11.75,-3.75) {0};

    \node at (1.25,-2.0) {$\cdots$};
    \node at (11.00,-2.0) {$\cdots$};
    % ... このように必要な座標に数字を配置していきます ...

\end{tikzpicture}
\caption{A frieze pattern of order $n$ arranged in staggered rows.}
\label{fig.1}
\end{figure}

%In their foundational work \cite{ConwayCoxeter1973}, Conway and Coxeter proved that any frieze pattern of order $n$, with all entries from row 2 to row $n-2$ being positive integers, corresponds uniquely to a triangulation of a convex $n$-gon, where the entries in the first non-trivial row of the frieze represent the number of triangles incident to each vertex of the polygon.
In their foundational work \cite{ConwayCoxeter1973},
Conway and Coxeter proved that every frieze pattern of order $n$
whose non-trivial entries $\{f_{i,i+k}\, ;\, i\in \mathbb{Z}, k=2,3,\ldots,n-2\}$
are positive integers corresponds uniquely to a triangulation of a convex
$n$-gon.
The entries in the first non-trivial row (``quiddity''\cite{ConwayCoxeter1973}) $\{f_{i,i+2}\, ;\, i=1,2,\ldots, n \}$
represent the number of triangles incident to each vertex of the polygon.
Since then, frieze patterns have been studied in connection with geometry
and combinatorics \cite{Genoud2015, Nisiyama2022}.
Recently, the introduction of cluster algebras \cite{FominZelevinsky2002}
brought increasing attention to the subject,
as Conway-Coxeter friezes give a natural realization of a cluster algebra of
type $A_{n-3}$.

As noted by Coxeter\cite{Coxeter1971}, if we assume that the entries $f_{i,j}$ are non-zero for
$i < j < i+n$,
the following property holds:
\begin{proposition}[cf.~\protect\cite{Coxeter1971}] \label{prop.coxeter}
For a frieze of order $n$, let $f_{i,j}$ be its entries.
If $f_{i,j} \neq 0$ for all $i < j < i+n$, then the Plücker relation
\begin{equation} \label{eq.plucker}
    f_{i,j}f_{k,\ell} - f_{i,k}f_{j,\ell} + f_{i,\ell}f_{j,k} = 0
\end{equation}
is satisfied for any integers $i \leq j \leq k \leq \ell \leq i+n$.
\end{proposition}
A proof of Proposition \ref{prop.coxeter} can be found in \cite{Knuth1997}
(answer to exercise 2.3.4.6-4).

The Plücker relations \eqref{eq.plucker} are fundamental to the geometry of
Grassmannians.
From this perspective, a frieze pattern can be viewed as a realization of the
coordinate ring of the Grassmannian
\begin{equation}
    \mathrm{Gr}(2,n) := GL_2(\mathbb{C})\backslash \left\{ Z=(v_1, v_2, \dots, v_n) \in \mathrm{Mat}_{2 \times n}(\mathbb{C}) \;\middle|\; \mathrm{rank}{Z}=2 \right\},
\end{equation}
where $v_i=(z_{1,i},z_{2,i})^T$ ($i=1, \ldots, n$).
In the next section, we introduce an identification 
of a frieze with a $2 \times n$
matrix $Z=(v_1,\dots,v_n)$ \cite{Nisiyama2022}.
This formulation enables us to analyze the symmetries of the frieze as
geometric transformations of the vector sequence
$\{v_i \mid i=1, \ldots, n\}$
induced by the action of the symmetric group $\mathfrak{S}_n$.
In this paper, we provide an explicit description of these induced
transformations.
%Furthermore, we investigate how the parity of $n$ dictates the algebraic constraints required to maintain the well-defined structure and internal consistency of the frieze.
If $n$ is odd, the entire symmetric group $\mathfrak{S}_n$ acts on the
frieze;
if $n$ is even, the dihedral group $\mathcal{D}_{n}$ ($\subset \mathfrak{S}_n$) acts on it.
%We provide an explicit description of the symmetry that emerges specifically in the case of order six, and we formulate a conjecture of the full symmetry group for the remaining even orders.
%We will show that in the case of odd order and even order frieze of odd-order admits S and

This paper is organized as follows.
Section \ref{section.2} establishes the algebraic foundations of frieze
patterns, describes their correspondence with the matrix space $\widehat{\mathcal{M}}_n$, and
introduces the normal forms.
Section 3 then investigates the induced action of the symmetric group
$\mathfrak{S}_n$.
For odd $n$, we demonstrate that the unique normal form yields a
well-defined action of $\mathfrak{S}_n$.
For even $n$, we analyze which subgroup can act.
%the constraints imposed by an additional multi-ratio invariant, provide an explicit description of the enhanced symmetry emerging specifically in order six, and formulate a conjecture for the full symmetry group in higher even orders.

%西山 (増田 大久保 津田 田内) 小木曽 井上 斉藤 野海 Willox for their  interest encouragement

\section{Frieze Patterns and the Matrix Space $\widehat{\mathcal{M}}_n$} \label{section.2}

In this section, we establish the algebraic foundations for friezes extended to
the entire planar lattice $\mathbb{Z} \times \mathbb{Z}$.

\begin{comment}
\begin{definition}[cf.~\protect\cite{Short2023}] \label{def.tame}
Let $R$ be a subset of an integral domain and $n \geq 4$ be an integer.
An array of numbers $F = (f_{i,j})_{i,j \in \mathbb{Z}}$ ($f_{i, j}\in R$) is called a frieze of order $n$ over $R$ if it satisfies the global Plücker relation,
\begin{equation} \label{eq:grobalPlucker}
f_{i,j}f_{k,\ell} - f_{i,k}f_{j,\ell} + f_{i,\ell}f_{j,k} = 0,
\end{equation}
and the boundary conditions
\begin{equation} \label{eq:boundary}
f_{i,i} = f_{i,i+n} = 0, \quad f_{i,i+1} = f_{i+1,i+n} = 1.
\end{equation}
In the following, we denote the set of all friezes of order $n$ over $\mathbb{C}$ by $\mathcal{CCF}_n$. We define $\mathcal{CCF}'_n$ as a subset of $\mathcal{CCF}_n$:
\begin{equation}
    \mathcal{CCF}'_n:=\{ F \in \mathcal{CCF}_{n} \mid f_{1,3} = 1 \}. \label{def:CCF'}
\end{equation}
\end{definition}
\end{comment}

\begin{definition}[cf.~\protect\cite{Short2023}] \label{def.tame}
Let $\mathbb{K}$ be a subset of the field of rational functions with complex
coefficients $\mathbb{C}(X)$, and let $n \ge 4$ be an integer.
An array of numbers $F = (f_{i,j})_{i,j \in \mathbb{Z}}$ ($f_{i, j}\in \mathbb{K}$) is called a frieze of order
$n$ over $\mathbb{K}$ if it satisfies the global Plücker relation,
\begin{equation} \label{eq:grobalPlucker}
f_{i,j}f_{k,\ell} - f_{i,k}f_{j,\ell} + f_{i,\ell}f_{j,k} = 0,
\end{equation}
and the boundary conditions,
\begin{equation} \label{eq:boundary}
f_{i,i} = f_{i,i+n} = 0, \quad f_{i,i+1} = f_{i+1,i+n} = 1.
\end{equation}
In what follows, we denote the set of all friezes of order $n$ over
$\mathbb{C}(X)$ by $\mathcal{CCF}_n$.
%We define $\mathcal{CCF}'_n$ as a subset of $\mathcal{CCF}_n$:
%\begin{equation}
    %\mathcal{CCF}'_n:=\{ F \in \mathcal{CCF}_{n} \mid f_{1,3} = 1 \}. \label{def:CCF'}
%\end{equation}
\end{definition}

\begin{remark}
With the initial conditions in \eqref{eq:boundary}, Definition
\ref{def.tame} is equivalent to the notion of a tame $\mathrm{SL}_2$-tiling
in \cite{Short2023}.
\end{remark}

Several noteworthy features of this definition deserve immediate attention,
particularly regarding how the global algebraic constraints recover and extend
classical properties of friezes.

By specializing the indices in the global Plücker relation
\eqref{eq:grobalPlucker} and applying the boundary conditions
\eqref{eq:boundary}, we obtain
\begin{align}
k = i+1, \, \ell = j+1 :\; &
f_{i,j} = -f_{j,i} \quad \text{(skew symmetry)}, \label{eq.skewSym}\\
k = i, \, \ell = i+1 :\; &
f_{i,j}f_{i+1,j+1} - f_{i+1,j}f_{i,j+1} = 1 \quad \text{(unimodular rule)}, \label{eq.unimodular} \\
k = i+1, \, \ell = i+n :\; &
f_{i,j} = -f_{i+n, j} \quad \text{(anti-periodicity in $i$)}, \label{eq.antiPeriodI} \\
k = j+1, \, \ell = j+n :\; &
f_{i,j} = -f_{i, j+n} \quad \text{(anti-periodicity in $j$)} . \label{eq.antiPeriodJ}
\end{align}

%First, by specializing the indices to $k = i$ and $\ell = i+1$ in the global Plücker relation \eqref{eq:grobalPlucker} and the applying the boundary conditions \eqref{eq:boundary}, we obtain the skew symmetry \begin{equation} \label{eq.skewSym}
%f_{i,j} = -f_{j,i}.
%\end{equation}
%Second, by setting $k = i+1$ and $\ell = j+1$, the global Plücker relation \eqref{eq:grobalPlucker} reduces to the unimodular rule (\ref{eq:diamond-determinant}).
%Third, the global Plücker relation \eqref{eq:grobalPlucker} inherently dictates the periodicity and symmetry of the frieze structure. By specializing the indices to $k = i+1$ and $\ell = i+n$, we obtain the anti-periodicity relation
%\begin{equation} \label{eq.antiPeriod}
    %f_{i,j} = -f_{i+n, j}.
%\end{equation}
%A similar specialization yields the counterpart for the second index, $f_{i,j} = -f_{i, j+n}$.
By combining the skew symmetry \eqref{eq.skewSym} with the anti-periodicity
\eqref{eq.antiPeriodI}, we arrive at the identity:
\begin{equation} \label{eq.glideSym.}
    f_{i,j} = f_{j, i+n}.
\end{equation}
This relation is called a \textit{glide symmetry}, a well-known classical
property of frieze patterns \cite{Genoud2015, Nisiyama2022}.
%Moreover, a significant advantage of this formal definition lies in its global consistency: it ensures that the frieze can be anti-periodically extended without contradiction across the entire planar lattice $\mathbb{Z} \times \mathbb{Z}$. This formulation allows us to treat the frieze as a global geometric object, where the twisted periodicity of $2n$ arises naturally from the underlying algebraic constraints of the integral domain.

Figure \ref{figure2} illustrates the anti-periodicity \eqref{eq.antiPeriodI}
and glide symmetry \eqref{eq.glideSym.} described above.
\begin{figure}[htbp]
\centering
\begin{tikzpicture}
    \draw[dashed] (0.12,-0.0) -- (3.88,0) -- (7.64,-2.5) -- (3.88,-2.5) -- cycle;
    \draw (-0.5,1.5) -- (-0.5,-3.5);
    \draw (-1.25,1.0) -- (9.0,1.0);

    \draw[dashed] (0.12,0) -- (0.12,-2.5);
    \draw[dashed] (3.88,0) -- (3.88,-2.5);
    \draw[dashed] (7.64,0) -- (7.64,-2.5);
    \draw[dashed] (3.88,0) -- (7.64,0);
    \draw[dashed] (0.12,-2.5) -- (3.88,-2.5);

    \node at (-0.0,1.25) {$\cdots$};
    \node at (0.5,1.25) {1};
    \node at (1.25,1.25) {2};
    \node at (2.0,1.25) {3};
    \node at (2.75,1.25) {4};
    \node at (3.5,1.25) {5};
    \node at (4.25,1.25) {6};
    \node at (5.0,1.25) {7};
    \node at (5.75,1.25) {8};
    \node at (6.5,1.25) {9};
    \node at (7.25,1.25) {10};
    \node at (8.0,1.25) {11};
    \node at (8.75,1.25) {$\cdots$};

    \node at (-1.0,0.5) {$\vdots$};
    \node at (-1.0,-0.25) {1};
    \node at (-1.0,-0.75) {2};
    \node at (-1.0,-1.25) {3};
    \node at (-1.0,-1.75) {4};
    \node at (-1.0,-2.25) {5};
    \node at (-1.0,-2.75) {6};
    \node at (-1.0,-3.25) {$\vdots$};

    \node at (0.5,-0.25) {0};
    \node at (1.25,-0.25) {1};
    \node at (2.0,-0.25) {2};
    \node at (2.75,-0.25) {1};
    \node at (3.5,-0.25) {1};
    \node at (4.25,-0.25) {0};
    \node at (5.0,-0.25) {-1};
    \node at (5.75,-0.25) {-2};
    \node at (6.5,-0.25) {-1};
    \node at (7.25,-0.25) {-1};
    \node at (8.0,-0.25) {0};

    \node at (0.5,-0.75) {-1};
    \node at (1.25,-0.75) {0};
    \node at (2.0,-0.75) {1};
    \node at (2.75,-0.75) {1};
    \node at (3.5,-0.75) {2};
    \node at (4.25,-0.75) {1};
    \node at (5.0,-0.75) {0};
    \node at (5.75,-0.75) {-1};
    \node at (6.5,-0.75) {-1};
    \node at (7.25,-0.75) {-2};
    \node at (8.0,-0.75) {-1};

    \node at (0.5,-1.25) {-2};
    \node at (1.25,-1.25) {-1};
    \node at (2.0,-1.25) {0};
    \node at (2.75,-1.25) {1};
    \node at (3.5,-1.25) {3};
    \node at (4.25,-1.25) {2};
    \node at (5.0,-1.25) {1};
    \node at (5.75,-1.25) {0};
    \node at (6.5,-1.25) {-1};
    \node at (7.25,-1.25) {-3};
    \node at (8.0,-1.25) {-2};

    \node at (0.5,-1.75) {-1};
    \node at (1.25,-1.75) {-1};
    \node at (2.0,-1.75) {-1};
    \node at (2.75,-1.75) {0};
    \node at (3.5,-1.75) {1};
    \node at (4.25,-1.75) {1};
    \node at (5.0,-1.75) {1};
    \node at (5.75,-1.75) {1};
    \node at (6.5,-1.75) {0};
    \node at (7.25,-1.75) {-1};
    \node at (8.0,-1.75) {-1};

    \node at (0.5,-2.25) {-1};
    \node at (1.25,-2.25) {-2};
    \node at (2.0,-2.25) {-3};
    \node at (2.75,-2.25) {-1};
    \node at (3.5,-2.25) {0};
    \node at (4.25,-2.25) {1};
    \node at (5.0,-2.25) {2};
    \node at (5.75,-2.25) {3};
    \node at (6.5,-2.25) {1};
    \node at (7.25,-2.25) {0};
    \node at (8.0,-2.25) {-1};

    \node at (0.5,-2.75) {0};
    \node at (1.25,-2.75) {-1};
    \node at (2.0,-2.75) {-2};
    \node at (2.75,-2.75) {-1};
    \node at (3.5,-2.75) {-1};
    \node at (4.25,-2.75) {0};
    \node at (5.0,-2.75) {1};
    \node at (5.75,-2.75) {2};
    \node at (6.5,-2.75) {1};
    \node at (7.25,-2.75) {1};
    \node at (8.0,-2.75) {0};

    \node at (2.0,0.5) {$\vdots$};
    \node at (2.0,-3.25) {$\vdots$};
    \node at (5.75,0.5) {$\vdots$};
    \node at (5.75,-3.25) {$\vdots$};
    \node at (-0.13,-1.25) {$\cdots$};
    \node at (8.75,-1.25) {$\cdots$};
    % ... このように必要な座標に数字を配置していきます ...

\end{tikzpicture}
\caption{The extended frieze pattern of order 5.}
\label{figure2}
\end{figure}

Henceforth, we let $\mathbb{K}$ be the field of rational functions with complex
coefficients $\mathbb{C}(X)$.

%By specializing the indices of the global Plücker relation \eqref{eq:grobalPlucker} to $k=1$ and $\ell=2$, and invoking the normalization \eqref{eq:boundary}, any entry $f_{i,j}$ of the frieze can be expressed as a $2 \times 2$ determinant:
%\begin{equation} \label{eq:matrix-representation}
    %f_{i,j} = \begin{vmatrix}
    %f_{1,i} & f_{1,j} \\
    %f_{2,i} & f_{2,j}
    %\end{vmatrix}.
%\end{equation}
 %This representation implies that the entire frieze is uniquely determined by its first two rows. Consequently, given a $2 \times n$ matrix whose columns are vectors in $\mathbb{C}^2$, the remaining entries of the frieze are obtained as the determinants of the corresponding $2 \times 2$ minors. To formalize this correspondence, we define the following sets of matrices.

Let $\mathcal{M}_n$ be the set of $2 \times n$ matrices such that all $2 \times 2$
minors are non-zero:
\begin{equation}
\mathcal{M}_n :=
\left\{
    Z=\begin{pmatrix}
        z_{1,1} & z_{1,2}  & \cdots & z_{1,n}\\
        z_{2,1} & z_{2,2}& \cdots & z_{2,n}
    \end{pmatrix}
    \;\middle|\; D_{i, j}(Z)\ne 0 \text{ for } 1\leq i, j\leq n
    \right\}, \label{def:M}
\end{equation}
%\begin{equation} \label{def:M}
    %\mathcal{M}_n := \left\{ (v_1, v_2, \dots, v_n) \in \mathrm{Mat}_{2 \times n}(\mathbb{C}) \;\middle|\; D_{i,j} \neq 0 \right\},
%\end{equation}
where
\begin{equation}
D_{i,j}(Z) := \det
\begin{pmatrix}
z_{1,i} & z_{1,j}\\
z_{2,i} & z_{2,j}
\end{pmatrix} \label{dfn:Dij}
\end{equation}
for $1 \leq i, j \leq n$.
We omit the argument $Z$ when there is no risk of confusion.
%To formalize the correspondence between friezes and matrices, we define a specific set of representative matrices.
For a matrix $Z\in\mathcal{M}_{n}$ as in \eqref{def:M}, we define the coordinate function
$x_{i,j}$ by
\begin{equation}
    x_{i,j}(Z):=z_{i,j}. \label{def:coordinateFunc}
\end{equation}
We define $\widehat{\mathcal{M}}_n$ and $\widehat{\widehat{\mathcal{M\,}}}_n$ as subsets of $\mathcal{M}_n$ as follows:
\begin{align}
    \widehat{\mathcal{M}}_n &:= \left\{
    \begin{pmatrix}
        0 & 1 & * & * & \cdots & * & 1 \\
        -1 & 0 & 1 & * & \cdots & * & *
    \end{pmatrix}\in
    \mathcal{M}_{n}
    \;\middle|\; D_{i,i+1} = 1 \text{ for } 1 \leq i \leq n, \text{ and }
    D_{1,n} = 1
    \right\}, \label{def:hatM}\\
    \widehat{\widehat{\mathcal{M\,}}}_{n} &:= \{ Z \in \widehat{\mathcal{M}}_n
    \mid z_{1,3} = 1 \}. \label{def:hatM'}
\end{align}
%An element $\widehat{Z}(x_1,\cdots,x_{n-3}):=()\in \widehat{\mathcal{M}}_n$ can be parametlize as follows: x_1,\cdots,x_{n-3}.
We remark that $\widehat{\widehat{\mathcal{M\,}}}_{n} \subset \widehat{\mathcal{M}}_n \subset \mathcal{M}_n$.
When $n=2m$ and $c\in\mathbb{C}^\times$, we define a map
$\eta(c): \mathcal{M}_{2m}  \to  \mathcal{M}_{2m}$,
\begin{equation}
\eta(c)(Z)=
\mathrm{diag}(c, 1)\cdot Z \cdot\mathrm{diag}(1, c^{-1},\ldots, 1, c^{-1}).
\end{equation}
For $Z\in \widehat{\mathcal{M}}_{2m}$, $z_{1,3}\ne 0$ because $D_{1,3}(Z)\ne 0$. Thus, $\eta(z_{1,3}^{-1})(Z)\in \widehat{\widehat{\mathcal{M\,}}}_{n}$. We define
a map $p: \widehat{\mathcal{M}}_{2m}\rightarrow \widehat{\widehat{\mathcal{M\,}}}_{2m}$ by
\begin{equation}
p(Z)=\eta(z_{1,3}^{-1})(Z).
\end{equation}
For $Z\in \widehat{\widehat{\mathcal{M\,}}}_{2m}$, it is clear that $p \circ \eta(c)(Z)=Z$ ($c\in\mathbb{K}^\times$).

We now introduce a map $\alpha$ as follows: %$ : \mathcal{CCF}_n \to \widehat{\mathcal{M}}_n$ and $\beta : \widehat{\mathcal{M}}_n \to \mathcal{CCF}_n$ by
\begin{equation}
\begin{array}{ccccc}
\alpha & : & \mathcal{CCF}_n & \to & \widehat{\mathcal{M}}_n \\
&& (f_{i,j})_{i,j \in \mathbb{Z}} & \mapsto & (f_{i,j})_{i=1,2, \, j=1,2,\dots,n}.
\end{array}
\end{equation}

\begin{proposition}
The inverse map $\alpha^{-1}$ can be defined as follows:
\begin{equation}
(z_{i, j})_{i=1, 2, j=1, 2,...,n} \mapsto (f_{i,j})_{i,j \in \mathbb{Z}},
\end{equation}
where
\begin{equation} \label{eq:alphaInv}
f_{i, j} = (-1)^{k+\ell}(z_{1, i'}z_{2, j'}-z_{1, j'}z_{2, i'}), \quad i=kn+i', \quad j=\ell n+j', \quad 1 \leq i', j' \leq n.
\end{equation}
\end{proposition}

\begin{proof}
For $i=1$ and $1 \leq j \leq n$, since $k=\ell=0$, formula \eqref{eq:alphaInv}
reduces to
\begin{equation}
f_{1,j}=z_{1,1}z_{2,j}-z_{1,j}z_{2,1}=z_{1,j}.
\end{equation}
Similarly, for $i=2$ and $1 \leq j \leq n$, we obtain $f_{2, j}=z_{2, j}$.
Next, we consider the case $1 \leq i, j \leq n$. By specializing the indices of the
global Plücker relation \eqref{eq:grobalPlucker} to $k=1$ and
$\ell=2$ and invoking the normalization \eqref{eq:boundary}, we can express
$f_{i,j}$ as a $2 \times 2$ determinant:
\begin{equation} \label{eq:matrix-representation}
    f_{i,j} = \begin{vmatrix}
    f_{1,i} & f_{1,j} \\
    f_{2,i} & f_{2,j}
    \end{vmatrix}=
    \begin{vmatrix}
    z_{1,i} & z_{1,j} \\
    z_{2,i} & z_{2,j}
    \end{vmatrix}.
\end{equation}
The formula \eqref{eq:alphaInv} follows from the anti-periodicity relations
\eqref{eq.antiPeriodI} and \eqref{eq.antiPeriodJ}.
\end{proof}

To relate $\mathcal{M}_n$ and $\widehat{\mathcal{M}}_n$, we consider the quotient space
\begin{equation}
GL_2(\mathbb{K})\backslash \mathcal{M}_n/(\mathbb{K}^\times)^n \label{eq:quotient}
\end{equation}
by the equivalence relation defined as follows: for $Z_1, Z_2 \in \mathcal{M}_n$, we say
$Z_1 \sim Z_2$ if and only if $Z_1=LZ_2 R$ for $L\in GL_2(\mathbb{K})$ and $R=\mathrm{diag}(d_1, \ldots, d_n)$
($d_1, \ldots, d_n \in \mathbb{K}^\times$).
We remark that this quotient space coincides with that used in the study of
hypergeometric functions \cite{GGR1992,Kimura2007}.

When $n$ is even, for any matrix $Z \in \mathcal{M}_{2m}$, there exists an invariant
$\mathrm{mr}(Z)$ with respect to the action of $GL_2(\mathbb{K}) \times (\mathbb{K}^\times)^{2m}$:
\begin{equation}
    \mathrm{mr}(Z) := \frac{D_{2,3} D_{4,5} \cdots D_{2m-2,2m-1} D_{1,2m}}{D_{1,2} D_{3,4} \cdots D_{2m-3,2m-2}D_{2m-1,2m}}. \label{eq:Phi}
\end{equation}
Using this invariant, we define a subset $\mathcal{M}'_{2m}$ of $\mathcal{M}_{2m}$ as follows:
\begin{equation} \label{def:M'}
    \mathcal{M}'_{2m} := \{ Z \in \mathcal{M}_{2m} \mid \mathrm{mr}(Z) = 1 \}.
\end{equation}

The following lemma is an immediate consequence of \eqref{def:hatM}.
\begin{lemma}
For every $Z \in \widehat{\mathcal{M}}_{2m}$, we have $\mathrm{mr}(Z)=1$.
\end{lemma}
We remark that $\widehat{\widehat{\mathcal{M\,}}}_{2m} \subset \widehat{\mathcal{M}}_{2m} \subset \mathcal{M}'_{2m}$.

%By the definitions \eqref{def:M}, \eqref{def:hatM}, \eqref{def:hatM'} and \eqref{def:M'}, it holds $\widehat{\mathcal{L}}_{n} \subset \widehat{\mathcal{M}}_{n} \subset \mathcal{M}'_{n} \subset \mathcal{M}_{n}$.
%By the definitions \eqref{def:M} and \eqref{def:hatM'}, it holds $\mathcal{L}_{n} \subset \mathcal{M}_{n}$.

\begin{remark}
When we specialize $\mathbb{K}$ to $\mathbb{C}$, the quotient space
\eqref{eq:quotient} can be interpreted in two different ways, depending on the
order in which the quotients are taken.
On the one hand, taking the quotient by $GL_2(\mathbb{C})$ first yields the quotient
of the Grassmannian $\mathrm{Gr}(2, n)$ by the right action of $(\mathbb{C}^\times)^n$.
On the other hand, taking the quotient by $(\mathbb{C}^\times)^n$ first yields the quotient
of the configuration space of $n$ points on the complex projective line
$\mathbb{P}^1$ by the left action of $GL_2(\mathbb{C})$
(cf. \cite{Aomoto1977}).
The quotient space $\mathcal{M}'_{2m}/(\mathbb{C}^\times)^{2m}$ was introduced in \cite{Genoud2021}, where it is
denoted by $\mathcal{C}^{+}_{2m}$.
The invariant $\mathrm{mr}(Z)$ in \eqref{eq:Phi} is known as the multi-ratio
associated with a configuration of $2m$ points on the projective line
$\mathbb{P}^1$.
In the six-point case, this multi-ratio was studied by Bobenko et al. in their
work on integrable circle patterns \cite{BHS2002,BobenkoHoffmann2001}.
\end{remark}

%\begin{thm} \label{prop:normal-form}
%Let $n \geq 5$ be an odd integer. For any $Z \in \mathcal{M}_n$, there exists a unique matrix $\widehat{Z} \in \widehat{\mathcal{M}}_n$ such that $Z \sim \widehat{Z}$. We refer to $\widehat{Z}$ as the normal form of $Z$.
%\end{thm}

\begin{thm} \label{thm:normalize}
    Let $n \geq 4$ be an integer and $Z \in \mathcal{M}_n$.
\begin{enumerate}[label=\textup{(\roman*)}]
    \item For odd $n$, there exists a unique matrix $\widehat{Z} \in \widehat{\mathcal{M}}_n$ such
    that $Z \sim \widehat{Z}$.
    %\item For even $n$, there exists a unique matrix $\widehat{Z} \in \widehat{\widehat{\mathcal{M\,}}}_{n}$ such that $Z \sim \widehat{Z}$ if and only if $Z \in \mathcal{M}'_n$.
    \item For even $n$, a unique matrix $\widehat{Z} \in \widehat{\widehat{\mathcal{M\,}}}_{n}$ such that
    $Z \sim \widehat{Z}$ exists if and only if $Z \in \mathcal{M}'_n$.
\end{enumerate}
We refer to this $\widehat{Z}$ as the normal form of $Z$.
\end{thm}

\begin{proof}
We fix a matrix $Z = (v_1, v_2, \dots, v_n) \in \mathcal{M}_n$, where $v_i = (a_i, b_i)^\top$
($i=1,\ldots, n$). Let $L_1, L_2$ and $R$ be given by
\begin{equation}
\begin{aligned}
    &L_1 = \begin{pmatrix} -b_1 & a_1 \\ -b_2 & a_2 \end{pmatrix}, \quad
    L_2 = \mathrm{diag}(\mu_1, \mu_2), \quad
    R = \mathrm{diag}(\lambda_1, \dots, \lambda_n), \\
&\mu_1, \mu_2, \lambda_1, \ldots, \lambda_n \in \mathbb{K}^\times,
\end{aligned}
\end{equation}
Then we have
%\begin{equation}
%L_1 Z=
%\begin{pmatrix}
% 0 & D_{1,2} & D_{1,3} & \cdots & D_{1,n} \\
% -D_{1,2} & 0 & \mu_2 \lambda_3 D_{2,3} & \cdots & \mu_2 \lambda_n D_{2,n} \end{pmatrix}.
%\end{equation}
%Then, $L_1$ is an element of $GL_2(\mathbb{C})$, and the columns of $L_1 Z$ are expressed as $(D_{1,j}(Z), D_{2,j}(Z))^\top$. By setting $L_2 = \mathrm{diag}(\mu_1, \mu_2) $ and $R = \mathrm{diag}(\lambda_1, \dots, \lambda_n) \in (\mathbb{C}^\times)^n$, the resulting matrix $\widehat{Z} = L_2 L_1 Z R$ takes the form:
\begin{equation}
    L_2 L_1 ZR = \begin{pmatrix} 0 & \mu_1 \lambda_2 D_{1,2} & \mu_1 \lambda_3 D_{1,3} & \cdots & \mu_1 \lambda_n D_{1,n} \\ -\mu_2 \lambda_1 D_{1,2} & 0 & \mu_2 \lambda_3 D_{2,3} & \cdots & \mu_2 \lambda_n D_{2,n} \end{pmatrix}, \label{eq:zTilde}
\end{equation}
where $D_{i,j} = D_{i,j}(Z)$.
For simplicity, we write $\widetilde{Z}:=L_2 L_1 ZR$.
To satisfy the normalization conditions \eqref{def:hatM} for the first three
columns of $\widehat{\mathcal{M}}_n$, the values of $\lambda_1, \lambda_2$, and
$\lambda_3$ are uniquely determined as
$\lambda_1 = 1/(\mu_2 D_{1,2})$, $\lambda_2 = 1/(\mu_1 D_{1,2})$, and $\lambda_3 = 1/(\mu_2 D_{2,3})$.
Furthermore, the conditions \eqref{def:hatM} for adjacent minors,
\begin{equation}
\mu_1 \mu_2 \lambda_i \lambda_{i+1} (D_{1,i}D_{2,i+1}-D_{1,i+1}D_{2,i}) = \mu_1 \mu_2 \lambda_i \lambda_{i+1} D_{1,2}D_{i,i+1} = 1
\end{equation}
allow us to recursively solve for $\lambda_k$. For $k = 4, \dots, n$, the general
form of $\lambda_k$ is given by:
\begin{equation}
    \lambda_k =
    \begin{cases}
    \displaystyle \frac{1}{\mu_2} \cdot \frac{D_{3,4}D_{5,6} \cdots D_{k-2,k-1}}{D_{2,3} D_{4,5} \cdots D_{k-1,k}} & (\text{if } k \text{ is odd}), \\[15pt]
    \displaystyle \frac{1}{\mu_1} \cdot \frac{D_{2,3} D_{4,5} \cdots D_{k-2,k-1}}{D_{1,2} D_{3,4} \cdots D_{k-1,k}} & (\text{if } k \text{ is even}).
    \end{cases}
\end{equation}

%First, we consider case (i). Let $n$ be an odd integer.
By substituting these values into \eqref{eq:zTilde}, the matrix $\widetilde{Z}$ is
explicitly written as:
\begin{equation}
\widetilde{Z} =
\begin{cases}
\begin{pmatrix}
0 & 1 & \displaystyle \frac{\mu_1 D_{1,3}}{\mu_2 D_{2,3}} & \displaystyle \frac{D_{2,3}D_{1,4}}{D_{1,2}D_{3,4}} & \cdots & \displaystyle \frac{\mu_1 D_{3,4}D_{5,6}\cdots D_{n-2,n-1}D_{1,n}}{\mu_2 D_{2,3}D_{4,5}\cdots D_{n-3,n-2}D_{n-1,n}} \\[12pt]
-1 & 0 & 1 & \displaystyle \frac{\mu_2 D_{2,3}D_{2,4}}{\mu_1 D_{1,2}D_{3,4}} & \cdots & \displaystyle \frac{D_{3,4}D_{5,6}\cdots D_{n-2,n-1}D_{2,n}}{D_{2,3}D_{4,5}\cdots D_{n-3,n-2}D_{n-1,n}}
\end{pmatrix} & (\text{if } n \text{ is odd}), \\[30pt]
\begin{pmatrix}
0 & 1 & \displaystyle \frac{\mu_1 D_{1,3}}{\mu_2 D_{2,3}} & \displaystyle \frac{D_{2,3} D_{1,4}}{D_{1,2} D_{3,4}} & \cdots & \mathrm{mr}(Z) \ \\[15pt]
-1 & 0 & 1 & \displaystyle \frac{\mu_2 D_{2,3} D_{2,4}}{\mu_1 D_{1,2} D_{3,4}} & \cdots &  \displaystyle \frac{\mu_2 D_{2,3}D_{4,5}\cdots D_{n-2,n-1}D_{2,n}}{\mu_1 D_{1,2}D_{3,4}\cdots D_{n-3,n-2}D_{n-1,n}}
\end{pmatrix} & (\text{if } n \text{ is even}),
\end{cases} \label{eq:hatZ}
\end{equation}
where both cases satisfy the conditions in \eqref{def:hatM} except for the
$(1,n)$-entry.
If $n$ is odd, we can set
%From the condition that the $(1,n)$-entry of $\widehat{Z}$ is $1$, we can determine the value of $\mu_1/\mu_2$ as
\begin{equation}
    \frac{\mu_1}{\mu_2} = \frac{D_{2,3} D_{4,5} \cdots D_{n-3,n-2}D_{n-1,n}}{D_{3,4} D_{5,6} \cdots D_{n-2,n-1} D_{1,n}}
\end{equation}
to make the $(1,n)$-entry of $\widetilde{Z}$ equal to $1$.
If $n$ is even, we require the condition $\mathrm{mr}(Z)=1$, which means that
$Z \in \mathcal{M}'_n$. Furthermore, we set
\begin{equation}
    \frac{\mu_1}{\mu_2} = \frac{D_{2,3}}{D_{1,3}}
\end{equation}
to make the $(1,3)$-entry of $\widetilde{Z}$ equal to $1$, which means
$\widetilde{Z} \in \widehat{\widehat{\mathcal{M\,}}}_{n}$.

%By setting $\mu_1 = D_{2,3} D_{4,5} \cdots D_{n-3,n-2}D_{n-1,n}$ and $\mu_2 = D_{3,4} D_{5,6} \cdots D_{n-2,n-1} D_{1,n}$, it follows that $L_2 \in GL_2(\mathbb{C})$ and $\widehat{Z} \in \widehat{\mathcal{M}}_n$. This completes the construction of $\widehat{Z}$.

We then address the uniqueness of the normal form. Suppose there exist two
representatives $\widehat{Z}_1, \widehat{Z}_2 \in \widehat{\mathcal{M}}_n$ (resp. $\widehat{\widehat{\mathcal{M\,}}}_{n}$) such that $Z \sim \widehat{Z}_1 \sim \widehat{Z}_2$. This
implies $\widehat{Z}_2 = L \widehat{Z}_1 R$ for some $L \in GL_2(\mathbb{K})$ and $R = \mathrm{diag}(\lambda_1, \dots, \lambda_n)$.
From the first and second columns of \eqref{def:hatM}, we see that $L$
must be a diagonal matrix $\mathrm{diag}(1/\lambda_2, 1/\lambda_1)$. Since $\widehat{Z}_2 \in \widehat{\mathcal{M}}_n$, we obtain the
following $n$ relations:
\begin{align}
1&=D_{i,i+1}(\widehat{Z}_2) = (\lambda_i \lambda_{i+1}/\lambda_1 \lambda_2) D_{i,i+1}(\widehat{Z}_1) = \lambda_i \lambda_{i+1}/\lambda_1 \lambda_2 \quad(1\leq i\leq n-1), \label{eq:Di,i+1}\\
1&=D_{1,n}(\widehat{Z}_2)=(\lambda_1 \lambda_{n}/\lambda_1 \lambda_2) D_{1,n}(\widehat{Z}_1) = \lambda_n/\lambda_2. \label{eq:D1,n}
\end{align}
Moreover, if $\widehat{Z}_2 \in \widehat{\widehat{\mathcal{M\,}}}_{n}$, we obtain the following relation from the condition
on the $(1,3)$-entry:
\begin{equation}
1=\frac{\lambda_3}{\lambda_2}. \label{eq:(1,3)-entry}
\end{equation}
The equations \eqref{eq:Di,i+1} imply the recurrence $\lambda_{i+1} = \lambda_1 \lambda_2 / \lambda_i$, which yields
\begin{equation}
\lambda_1=\lambda_3=\lambda_5=\cdots, \quad \lambda_2=\lambda_4=\lambda_6=\cdots. \label{eq:lambda}
\end{equation}
If $n$ is odd, it follows from \eqref{eq:lambda} that $\lambda_n = \lambda_1$.
Then, by \eqref{eq:D1,n}, we have $\lambda_1 = \lambda_2 = \dots = \lambda_n$.
On the other hand, if $\widehat{Z}_2 \in \widehat{\widehat{\mathcal{M\,}}}_{n}$, then \eqref{eq:(1,3)-entry} yields
$\lambda_2 = \lambda_3$, and again by \eqref{eq:lambda}, we obtain $\lambda_1 = \lambda_2 = \dots = \lambda_n$. Therefore,
in either case, $L$ and $R$ can be expressed as $L = (1/\lambda) I$ and
$R = \lambda I$ for some $\lambda \in \mathbb{K}^\times$, which implies $\widehat{Z}_2 = L \widehat{Z}_1 R = \widehat{Z}_1$.
\end{proof}

Let $\varphi_{\text{odd}}: \mathcal{M}_{2m+1} \rightarrow \widehat{\mathcal{M}}_{2m+1}$ (resp. $\varphi_{\text{even}}: \mathcal{M}'_{2m} \rightarrow \widehat{\mathcal{M}}_{2m}$) be the map that assigns $Z \in \mathcal{M}_{2m+1}$ (resp.
$\mathcal{M}_{2m}'$) to its normal form in the sense of Theorem \ref{thm:normalize}.
More explicitly, for $Z\in \mathcal{M}_{2m+1}$, we can write
\begin{equation}
\begin{aligned}
&\varphi_{\text{odd}}(Z):=
\begin{pmatrix}
0 & 1 & \displaystyle \widehat{z}_{1,3} & \widehat{z}_{1,4} & \cdots & \displaystyle \widehat{z}_{1,2m} & 1 \\[0pt]
-1 & 0 & 1 & \widehat{z}_{2,4} & \cdots & \widehat{z}_{2,2m} & \widehat{z}_{2,2m+1}
\end{pmatrix}, \\[15pt]
&\widehat{z}_{1,k} =
\begin{cases}
 \displaystyle\frac{D_{k+1,k+2} D_{k+3,k+4} \cdots D_{2m,2m+1} D_{1,k}}{D_{k,k+1} D_{k+2,k+3} \cdots D_{2m-1,2m} D_{1,2m+1}} & (\text{if } k \text{ is odd}),  \\[15pt]
\displaystyle\frac{D_{2,3} D_{4,5} \cdots D_{k-2,k-1} D_{1,k}}{D_{1,2} D_{3,4} \cdots D_{k-3, k-2}D_{k-1,k}} & (\text{if } k \text{ is even}),
\end{cases} \\[15pt]
&\widehat{z}_{2,k} =
\begin{cases}
\displaystyle\frac{D_{3,4} D_{5,6} \cdots D_{k-2,k-1} D_{2,k}}{D_{2,3} D_{4,5} \cdots D_{k-3, k-2}D_{k-1,k}} & (\text{if } k \text{ is odd}),  \\[15pt]
\displaystyle\frac{D_{k+1,k+2} D_{k+3,k+4} \cdots D_{2m-1,2m} D_{1,2m+1} D_{2,k}}{D_{k,k+1} D_{k+2,k+3} \cdots D_{2m-2,2m-2}D_{2m,2m+1} D_{1,2}} & (\text{if } k \text{ is even}),
\end{cases}
\end{aligned}\label{eq:hatZ1k}
\end{equation}
and for $Z'\in \mathcal{M}'_{2m}$,
\begin{equation}
\begin{aligned}
&\varphi_{\text{even}}(Z'):=\begin{pmatrix}
0 & 1 & 1 & \widehat{w}_{1,4} & \cdots & \widehat{w}_{1,2m-1} & 1 \ \\[0pt]
-1 & 0 & 1 & \widehat{w}_{2,4} & \cdots &  \widehat{w}_{2,2m-1} & \widehat{w}_{2,2m}
\end{pmatrix},
 \\[15pt]
&\widehat{w}_{1,k} =
\begin{cases}
 \displaystyle\frac{D_{3,4} D_{4,5} \cdots D_{k-2,k-1} D_{1,k}}{D_{1,3} D_{4,5} \cdots D_{k-3,k-2} D_{k-1,k}} & (\text{if } k \text{ is odd}),  \\[15pt]
\displaystyle\frac{D_{2,3} D_{4,5} \cdots D_{k-2,k-1} D_{1,k}}{D_{1,2} D_{3,4} \cdots D_{k-3, k-2}D_{k-1,k}} & (\text{if } k \text{ is even}),
\end{cases} \\[15pt]
&\widehat{w}_{2,k} =
\begin{cases}
\displaystyle\frac{D_{3,4} D_{5,6} \cdots D_{k-2,k-1} D_{2,k}}{D_{2,3} D_{4,5} \cdots D_{k-3, k-2}D_{k-1,k}} & (\text{if } k \text{ is odd}),  \\[15pt]
\displaystyle\frac{D_{1,3} D_{4,5} \cdots D_{k-2,k-1} D_{2,k}}{D_{1,2} D_{3,4} \cdots D_{k-3,k-2}D_{k-1,k}} & (\text{if } k \text{ is even}).
\end{cases}
\end{aligned}
\end{equation}

\begin{remark}
    A formula similar to the expressions for $\widehat{z}_{1,k}$ and $\widehat{z}_{2,k}$ in
    \eqref{eq:hatZ1k} has also appeared in the study of string theory (see
    \cite{HIS2012}, (2.10)).
\end{remark}

\section{Action of the symmetric group}
Based on Theorem \ref{thm:normalize}, the map $\varphi_{\text{odd}} : \mathcal{M}_{2m+1} \to \widehat{\mathcal{M}}_{2m+1}$, which assigns each
matrix $Z$ to its unique normal form $\widehat{Z}$, is well-defined.
For $\sigma\in \mathfrak{S}_{2m+1}$, define $\rho({\sigma}): \widehat{\mathcal{M}}_{2m+1} \to \widehat{\mathcal{M}}_{2m+1}$ as
\begin{equation}
    \rho(\sigma)(Z) = \varphi_{\text{odd}}(Z P_{\sigma}), \label{eq:action}
\end{equation}
where $P_{\sigma}$ is the permutation matrix corresponding to $\sigma$.
We define a map $\check{\rho}(\sigma): \mathcal{CCF}_{2m+1} \rightarrow \mathcal{CCF}_{2m+1}$ as
\begin{equation}
\check{\rho}(\sigma) := \alpha^{-1} \circ \rho(\sigma) \circ \alpha.
\end{equation}
This action yields a representation of the symmetric group on $\mathcal{CCF}_{2m+1}$.

In what follows, let $\sigma_{ij}$ denote the transposition $(i, j)$ for
$i, j = 1, \dots, 2m+1$.

\begin{example}
%The action of $\mathfrak{S}_5$ on $\mathcal{CCF}_{5}$ can be expressed by the matrices as $5 \times 5$ anti-symmetric matrices:
As an example, we consider the following $F\in\mathcal{CCF}_5$:
\begin{equation*}
\begin{array}{c@{}ccc@{}c}
& \vdots & \vdots&&\\
F\;=\;\cdots &
\begin{array}{ccccc}
            0 & 1 & x & y & 1 \\
            -1 & 0 & 1 & \frac{1+y}{x} & \frac{1+x+y}{xy} \\
            -x & -1 & 0 & 1 & \frac{1+x}{y} \\
            -y & -\frac{1+y}{x} & -1 & 0 & 1 \\
            -1 & -\frac{1+x+y}{xy} & -\frac{1+x}{y} & -1 & 0
\end{array}&
\begin{array}{ccccc}
            0 & -1 & -x & -y & -1 \\
            1 & 0 & -1 & -\frac{1+y}{x} & -\frac{1+x+y}{xy} \\
            x & 1 & 0 & -1 & -\frac{1+x}{y} \\
            y & \frac{1+y}{x} & 1 & 0 & -1 \\
            1 & \frac{1+x+y}{xy} & \frac{1+x}{y} & 1 & 0
\end{array}&\cdots\\
&\vdots&\vdots&
\end{array}
\end{equation*}
For simplicity, only the upper triangular part of the first $5 \times 5$ block
of $\check{\rho}(\sigma).F$ is shown below:
\begin{align*}
    \check{\rho}(\sigma_{12}).F &= \begin{pmatrix}
            0 & 1 & \frac{xy}{1+x+y} & -(1+y) & 1 \\
            & 0 & 1 & \frac{-(1+x+y)}{x} & \frac{1}{x} \\
            && 0 & 1 & \frac{-(1+x)}{1+x+y} \\
            &&& 0 & 1 \\
            &&&& 0
        \end{pmatrix}, \\
    \check{\rho}(\sigma_{23}).F &= \begin{pmatrix}
            0 & 1 & \frac{x}{1+y} & \frac{-y}{1+y} & 1 \\
            & 0 & 1 & \frac{1}{x} & \frac{-(1+x)(1+y)}{xy} \\
            && 0 & 1 & \frac{-(1+x+y)}{y} \\
            &&& 0 & 1 \\
            &&&& 0
        \end{pmatrix}, \\
    \check{\rho}(\sigma_{34}).F &= \begin{pmatrix}
            0 & 1 & -(1+x) & -(1+y) & 1 \\
            & 0 & 1 & \frac{y}{1+x} & \frac{-(1+x+y)}{(1+x)(1+y)} \\
            && 0 & 1 & \frac{x}{1+y} \\
            &&& 0 & 1 \\
            &&&& 0
        \end{pmatrix}, \\
    \check{\rho}(\sigma_{45}).F &= \begin{pmatrix}
            0 & 1 & \frac{-x}{1+x} & \frac{y}{1+x} & 1 \\
            & 0 & 1 & \frac{-(1+x+y)}{x} & \frac{-(1+x)(1+y)}{xy} \\
            && 0 & 1 & \frac{1}{y} \\
            &&& 0 & 1 \\
            &&&& 0
        \end{pmatrix}.
\end{align*}
\end{example}

In the case of $\varphi_{\text{even}}: \mathcal{M}'_{2m} \rightarrow \widehat{\widehat{\mathcal{M}\,}}_{2m}$, the action \eqref{eq:action} does not work in
general;
even if $\mathrm{mr}(Z')=1$ (i.e. $Z'\in\mathcal{M}'_{2m}$), the condition $\mathrm{mr}(Z'P_{\sigma})=1$ need not hold
for every $\sigma\in \mathfrak{S}_{2m}$.

Define $r$ and $t$ as follows:
\begin{equation}
r=\sigma_{1,2m}\sigma_{1,2m-1}\cdots\sigma_{12}, \quad t=\sigma_{2,2m}\sigma_{3,2m-1}\cdots\sigma_{m,m+2}.
\end{equation}
They generate the dihedral group $\mathcal{D}_{2m} \subset \mathfrak{S}_{2m}$.
By the definition of the multi-ratio \eqref{eq:Phi}, if $\mathrm{mr}(Z')=1$, then
\begin{equation}
    \mathrm{mr}(Z'P_{r})=\mathrm{mr}(Z'P_{t})=1.
\end{equation}
Thus we can define the action of $r, t$ on $\widehat{\widehat{\mathcal{M}\,}}_{2m}$:
\begin{equation}
    \rho(r)(Z') = \varphi_{\text{even}}(Z' P_{r}), \quad
    \rho(t)(Z') = \varphi_{\text{even}}(Z' P_{t}), \quad
    Z'\in \widehat{\widehat{\mathcal{M}\,}}_{2m}. \label{eq:dihedral}
\end{equation}

\begin{comment}
By direct calculation, we can show that the maps $\rho(r)$, $\rho(t)$, and $\eta(c)$ satisfy the following commutation relations:
\begin{equation}
\eta(c^{-1})\rho(r)=\rho(r)\eta(c), \quad
\eta(c)\rho(t)=\rho(t)\eta(c). \label{eq:swap}
\end{equation}
\end{comment}

For an element $\sigma$ of the dihedral group $\mathcal{D}_{2m}$ and a parameter
$c\in \mathbb{K}^\times$, we define the map $\rho(\sigma,c) \colon \widehat{\mathcal{M}}_{2m} \to \widehat{\mathcal{M}}_{2m}$ as follows:
\begin{equation}
    \rho(\sigma,c):=\eta(c)\circ\rho(\sigma)\circ p. \label{def:rhoEven}
\end{equation}

\begin{comment}
\begin{proposition}
For any parameters $c_1, \dots, c_{2m}\in\mathbb{K}^\times$ and $Z\in\widehat{\mathcal{M}}_{2m}$, the following relations hold:
\begin{align}
&\rho(r,c_{2m})\cdots\rho(r,c_{1})=
\eta\biggl(\frac{c_2\cdots c_{2m}}{c_1\cdots c_{2m-1}}\biggr), \\
&\rho(t,c_{2})\rho(t,c_{1})=
\eta(c_1 c_2), \\
&\rho(r,c_{4})\rho(t,c_{3})\rho(r,c_{2})\rho(t,c_{1})=
\eta\biggl(\frac{c_{1} c_{4}}{c_{2} c_{3}}\biggr).
\end{align}
\end{proposition}
\end{comment}

\begin{proposition}
For any parameters $c_1, \dots, c_{2m}\in\mathbb{K}^\times$ and $Z\in\widehat{\mathcal{M}}_{2m}$, the following relations hold:
\begin{align}
&\rho(r,c_{2m})\cdots\rho(r,c_{1})(Z)=
\eta\biggl(\frac{c_{2m}}{x_{1,3}(Z)}\biggr)(Z), \\
&\rho(t,c_{2})\rho(t,c_{1})(Z)=
\eta\biggl(\frac{c_2}{x_{1,3}(Z)}\biggr)(Z), \\
&\rho(r,c_{4})\rho(t,c_{3})\rho(r,c_{2})\rho(t,c_{1})(Z)=
\eta\biggl(\frac{c_{4}}{x_{1,3}(Z)}\biggr)(Z).
\end{align}
\end{proposition}

Let $i \colon \widehat{\widehat{\mathcal{M}\,}}_{2m} \hookrightarrow \widehat{\mathcal{M}}_{2m}$ be the natural inclusion map defined by $i(Z) = Z$.
We define an equivalence relation on $\widehat{\mathcal{M}}_{2m}$ as follows: for $Z_1, Z_2 \in \widehat{\mathcal{M}}_{2m}$,
$Z_1 \sim_{\text{mat}} Z_2$ if and only if there exists $c\in \mathbb{K}^{\times}$ such that $\eta(c)(Z_1) = Z_2$.
We introduce the quotient space $\widehat{\widehat{\mathcal{M}\,}}_{2m}/\sim_{\text{mat}}$ 
and a projection map 
$\pi: \widehat{\mathcal{M}'}_{2m} \rightarrow \widehat{\widehat{\mathcal{M}\,}}_{2m}/\sim_{\text{mat}}$. 
The following diagram illustrates these maps:
\begin{equation}
\begin{tikzcd}
\widehat{\widehat{\mathcal{M}\,}}_{2m} \arrow[r,hookrightarrow,"i"] &
\widehat{\mathcal{M}}_{2m} \arrow[r,rightarrow,"\pi"] &
\widehat{\mathcal{M}}_{2m}/\sim_{\text{mat}}. 
\end{tikzcd}
\end{equation}
It follows that $\widehat{\mathcal{M}}_{2m}/\sim_{\text{mat}} \ \cong \widehat{\widehat{\mathcal{M}\,}}_{2m}$.
The action of $\rho(\sigma,c)$ on $\widehat{\mathcal{M}}_{2m}/\sim_{\text{mat}}$ is independent of $c\in\mathbb{K}^\times$.
It follows that $\rho(\sigma,c)$ induces a representation of $\mathcal{D}_{2m}$ on
$\widehat{\mathcal{M}}_{2m}/\sim_{\text{mat}}$.

\begin{comment}
The map $\rho(\sigma,c)$ defined in \eqref{def:rhoEven} induces a map on the quotient space. In what follows, by abuse of notation, we shall also denote the induced map by $\rho(\sigma,c)$.
By the definition of the equivalence relation $\sim_{\text{mat}}$, any scaling map $\eta(c)$ acts as the identity on $\widehat{\mathcal{M}'}_{2m}/\sim_{\text{mat}}$.
Thus, we can simply write $\rho(\sigma,c)$ as $\bar{\rho}(\sigma)$.

\begin{align}
    \bar{\rho}(r)^{2m} &= \mathrm{id}, \\
    \bar{\rho}(t)^{2} &= \mathrm{id}, \\
    (\bar{\rho}(r)\bar{\rho}(t))^{2} &= \mathrm{id}.
\end{align}

We also define a map on $\check{\rho}(\sigma,c): \mathcal{CCF}_{2m} & \to & \mathcal{CCF}_{2m}$, where
\begin{equation}
\check{\rho}(\sigma,c) := \alpha^{-1} \circ \rho(r,c) \circ \alpha.
\end{equation}
\end{comment}

For $m=3$, the set $\mathcal{CCF}_{6}$ exhibits enhanced symmetry associated with
the Weyl group $W(B_3)$ ($\simeq W(C_3)$):
\begin{equation}
    \mathcal{D}_6 \subset W(B_3) \subset \mathfrak{S}_6.
\end{equation}
We define
\begin{equation}
    g_1=\sigma_{14}, \quad g_2=\sigma_{16}\sigma_{34}, \quad g_3=\sigma_{26}\sigma_{35}=t, \label{dfn:g123}
\end{equation}
which satisfy
\begin{equation}
g_i^2=\textrm{id} \; (i=1,2,3), \quad (g_1g_3)^2=(g_1g_2)^3=(g_2g_3)^4=\textrm{id}, \quad r=g_2g_1g_2g_3g_2, \label{eq:C3rels.}
\end{equation}
and hence $g_1$, $g_2$, and $g_3$ generate $W(B_3)$.

\begin{proposition}[Affolter, George and Ramassamy\protect\cite{AGR2025}, Lemma
4.2] \label{prop.mr}\label{prop.AGR}
    Let $Z$ be a matrix in $\mathcal{M}'_{6}$. For every permutation
    $\sigma$ in $\mathfrak{S}_6$ such that
    \begin{equation}
        \sigma(i+3) \equiv \sigma(i)+3 \mod 6 \quad (1\leq i \leq 6), \label{eq:mod6}
    \end{equation}
    we have $\mathrm{mr}(ZP_{\sigma})=1$, i.e. $ZP_\sigma\in\mathcal{M}'_{6}$.
\end{proposition}

\begin{remark}
An equivalent discussion can also be found in Proposition 2.1 of \cite{BobenkoHoffmann2001}.
\end{remark}

Since $g_1, g_2$, and $g_3$ in \eqref{dfn:g123} satisfy the condition
\eqref{eq:mod6}, Proposition \ref{prop.AGR} implies that $\mathrm{mr}(ZP_g) = 1$ holds for
all $Z \in \mathcal{M}'_{2m}$ and $g \in W(B_3)$.
Therefore, for $m = 3$, the Weyl group $W(B_3)$, which is larger than the
dihedral group $\mathcal{D}_6$, acts on
$\widehat{\mathcal{M}}_{6}/\sim_{\text{mat}} \ \cong \widehat{\widehat{\mathcal{M}\,}}_{6}$.

\begin{conjecture} \label{conj.dihedral}
If $m \ne 3$, the symmetry group of $\widehat{\widehat{\mathcal{M}\,}}_{2m}$ coincides with the dihedral
group $\mathcal{D}_{2m}$.
\end{conjecture}
Using the computer algebra systems SageMath, Maxima, and GAP, we have verified
that Conjecture \ref{conj.dihedral} holds for $m=2, 4, 5$, and $6$.

%where
In what follows, we define the map $\check{\rho}(\sigma,c): \mathcal{CCF}_{2m}\to\mathcal{CCF}_{2m}$ by
\begin{equation}
    \check{\rho}(\sigma,c):=\alpha^{-1}\circ \rho(\sigma,c)\circ \alpha.
\end{equation}
The construction of $\check{\rho}$ can be summarized as follows:

\begin{center}
    \begin{tikzcd}[
      row sep=1.4em,
      column sep=1.8em,
      ampersand replacement=\&
    ]
      \mathcal{CCF}_{2m}
        \arrow[r, "\alpha"]
        \arrow[d, "{\check{\rho}(\sigma,c)}"']
      \&
      \widehat{\mathcal{M}}_{2m}
        \arrow[r, "p"]
        \arrow[d, "{\rho(\sigma,c)}"]
      \&
      \widehat{\widehat{\mathcal{M}\,}}_{2m}
        \arrow[d, "\rho(\sigma)"]
      \\
      \mathcal{CCF}_{2m}
      \&
      \widehat{\mathcal{M}}_{2m}
        \arrow[l, "\alpha^{-1}"']
      \&
      \widehat{\widehat{\mathcal{M}\,}}_{2m}
        \arrow[l, "\eta(c)"']
     \&.
    \end{tikzcd}
\end{center}

\begin{example}
As an example, we consider the following $F\in\mathcal{CCF}_6$:
\begin{align*}
\begin{array}{c@{}cc}
& \vdots & \\
F\;=\;\cdots &
\begin{array}{cccccc}
            0 & 1 & x & y & z & 1 \\
            -1& 0 & 1 & \frac{1+y}{x} & \frac{x+z+yz}{xy} & \frac{(x+z)(1+y)}{xyz} \\
            -x&-1& 0 & 1 & \frac{x+z}{y} & \frac{x+xy+z}{yz} \\
            -y&-\frac{1+y}{x}&-1& 0 & 1 & \frac{1+y}{z} \\
            -z&-\frac{x+z+yz}{xy}&-\frac{x+z}{y}&-1& 0 & 1\\
        -1&-\frac{(x+z)(1+y)}{xyz}&-\frac{x+z+xy}{yz}&-\frac{1+y}{z}&-1& 0
\end{array}&\cdots\\
&\vdots&
\end{array}
\end{align*}
For simplicity, only the upper triangular part of the first $6 \times 6$ block
of $\check{\rho}(\sigma).F$ is shown below:
\begin{align*}
    \check{\rho}(r,c).F &= \begin{pmatrix}
            0 & 1 & \frac{c(1+y)}{x} & \frac{x+z+yz}{xy} & \frac{c(x+z)(1+y)}{xyz} & 1 \\
            & 0 & 1 & \frac{x+z}{cy} & \frac{x+z+xy}{yz} & \frac{x}{c} \\
            && 0 & 1 & \frac{c(1+y)}{z} & y \\
            &&& 0 & 1 & \frac{z}{c} \\
            &&&& 0 & 1\\
        &&&&& 0
        \end{pmatrix}, \\
    \check{\rho}(g_3, c).F = \check{\rho}(t, c).F &= \begin{pmatrix}
            0 & 1 & -cz & y & -cx & 1 \\
            & 0 & 1 & -\frac{1+y}{cz} & \frac{x+z+xy}{yz} & -\frac{(x+z)(1+y)}{cxyz} \\
            && 0 & 1 & -\frac{c(x+z)}{y} & \frac{x+z+yz}{xy} \\
            &&& 0 & 1 & -\frac{1+y}{cx} \\
            &&&& 0 & 1\\
        &&&&& 0
        \end{pmatrix}, \\
    \check{\rho}(g_1,c).F &= \begin{pmatrix}
            0 & 1 & -c & -\frac{y}{1+y} & -\frac{cx}{z} & 1 \\
            & 0 & 1 & -\frac{1}{c(1+y)} & -\frac{x+z+yz}{yz} & -\frac{x+z}{cxy} \\
            && 0 & 1 & -\frac{(1+x)(1+y)}{xy} & -\frac{1+x+y}{x} \\
            &&& 0 & 1 & \frac{y}{1+x} \\
            &&&& 0 & 1\\
        &&&&& 0
        \end{pmatrix}, \\
    \check{\rho}(g_2,c).F &= \begin{pmatrix}
            0 & 1 & -\frac{cx}{z} & -\frac{x+z+xy}{x+z} & \frac{cxy}{(x+z)(1+y)} & 1 \\
            & 0 & 1 & \frac{yz}{c(x+z)} & -\frac{x+z+yz}{(x+z)(1+y)} & \frac{1}{c} \\
            && 0 & 1 & \frac{cx}{(1+y)z} & -\frac{x+z}{z} \\
            &&& 0 & 1 & -\frac{1+y}{c} \\
            &&&& 0 & 1\\
        &&&&& 0
        \end{pmatrix}.
\end{align*}
\end{example}

\begin{comment}
For a parameter $c \in \mathbb{K}^\times$, we introduce a map $\eta_c$,
\begin{equation}
\begin{array}{ccccc}
\eta_c & : & \mathcal{M}'_{2m} & \to & \mathcal{M}'_{2m} \\
&& (x_{ij})_{i=1,2,j=1,\ldots,2m} & \mapsto & (c^{\epsilon(i,j)}x_{ij})_{{i=1,2,j=1,\ldots,2m}},
\end{array}
\end{equation}
where
\begin{equation} \label{eq:eta}
\epsilon(i,j) =
\begin{cases}
1 & (\text{if both } i \text{ and } j \text{ are odd}),  \\[0pt]
-1 & (\text{if both } i \text{ and } j \text{ are even}), \\[0pt]
0 & (\text{otherwise}).
\end{cases} \\
\end{equation}

We remark that a matrix $\eta(c)(Z')$ can be written as
\begin{equation}
\eta_c(Z')=
\mathrm{diag}(c, 1)\cdot Z' \cdot\mathrm{diag}(1, c^{-1},\ldots, 1, c^{-1}).
\end{equation}
\end{comment}

\begin{remark}
It is a natural direction to extend our matrix realization approach to
$SL_k$-frieze patterns by adopting the framework introduced by Surmann
\cite{Surmann2025}. Representing $SL_k$-frieze patterns via $k \times n$
matrices (or elements in $\mathrm{Gr}(k,n)$) provides a foundation for analyzing
actions of $\mathfrak{S}_n$ and its subgroups. Generalizing our classification of
symmetry groups to the $SL_k$ case remains an interesting topic for future
research.
\end{remark}

\section*{Acknowledgments}
The authors acknowledge Yasuyuki Hatsuda, Masatoshi Noumi and Yoshihisa Saito
for discussions on normal forms of elements in $\widehat{\mathcal{M}}_n$. Thanks are also due
to Rei Inoue, Takeyoshi Kogiso and Kyo Nishiyama for providing information on
frieze patterns.

During the preparation of this work, the authors used OpenAI Prism to assist in
the preparation and editing of this manuscript, and Google Gemini to improve
the clarity and style of the English, as well as to double-check the
mathematical computations and ensure accuracy.
After using these tools, the authors reviewed and edited the content as needed
and take full responsibility for the content of the publication.
%
%H. Masuda, H. Okubo, T. Tsuda, H. Tanaka, K. Oguiso, R. Inoue, H. Saito, Y. Nohara, and R. Willox for their interest and encouragement.


\begin{thebibliography}{9}

\bibitem{AGR2025}
N.C. Affolter, T. George and S. Ramassamy,
Integrable dynamics in projective geometry via dimers and triple crossing
diagram maps on the cylinder,
\textit{SIGMA} \textbf{21} (2025), 040.

\bibitem{Aomoto1977}
K. Aomoto, On the structure of integrals
of power product linear functions,
\textit{Sci. Papers. Gen. Univ. Tokyo} \textbf{27} (1977), 49--61.


\bibitem{BobenkoHoffmann2001}
A.I. Bobenko and T. Hoffmann,
Conformally symmetric circle packings: a generalization of Doyle's spirals,
\emph{Experiment. Math.} \textbf{10} (2001), 141--150.


\bibitem{BHS2002}
A.I. Bobenko, T. Hoffmann, and Yu.B. Suris, Hexagonal circle patterns and
integrable systems: Patterns with the multi-ratio property and Lax equations on
the regular triangular lattice, \textit{Internat. Math. Res. Notices}
\textbf{2002}, 111--164.

\bibitem{ConwayCoxeter1973}
J.H. Conway and H.S.M. Coxeter, Triangulated polygons and frieze patterns,
\emph{Math. Gaz.} \textbf{57} (1973), 87--94.

\bibitem{Coxeter1971}
H.S.M. Coxeter, Frieze patterns, \emph{Acta Arithmetica} \textbf{18} (1971),
297--310.

\bibitem{FominZelevinsky2002}
S. Fomin and A. Zelevinsky, Cluster algebras I: foundations, \emph{J. Amer.
Math. Soc.} \textbf{15} (2002), 497--529.

\bibitem{GGR1992}
I.M. Gel'fand, M.I. Graev and V.S. Retakh,
General hypergeometric systems of equations and series of hypergeometric type,
\textit{Russian Math. Surveys}, 47:4 (1992), 1--88.

\bibitem{HIS2012}
Y. Hatsuda, K. Ito and Y. Satoh,
T-functions and multi-gluon scattering amplitudes,
\textit{J. High Energy Phys.} \textbf{1202} (2012) 003.

\bibitem{Kimura2007}
H. Kimura, 
\textit{Introduction to Hypergeometric Functions}, 
Saieinsu-sha, 2007 [in Japanese]. 

\bibitem{Knuth1997}
D.E. Knuth,
\textit{The Art of computer programming, Vol. 1: fundamental Algorithms, Third
Edition}, Addison-Wesley Professional, 1997.

\bibitem{Genoud2015}
S. Morier-Genoud,
Coxeter’s frieze patterns at the crossroads of algebra, geometry and
combinatorics, \textit{Bull. London Math. Soc.} \textbf{47} (2015), 895--938.

\bibitem{Genoud2021}
S. Morier-Genoud,
Counting Coxeter’s friezes over a ﬁnite ﬁeld via moduli spaces,
\emph{Algebr. Comb.} \textbf{4} (2021), 225--240.

\bibitem{Nisiyama2022}
K. Nishiyama, \textit{Mathematics of frieze patterns --- A glimpse of geometry
and algebra}, Kyoritsu Shuppan Co., Ltd., 2022 [in Japanese].

\bibitem{Short2023}
I. Short, Classifying $\mathrm{SL}_2$-tilings, \emph{Trans. Amer. Math. Soc.}
\textbf{376} (2022), 1--38.

\bibitem{Surmann2025}
L. Surmann,
Counting tame $SL_3$- and $SL_4$-frieze patterns over finite fields,
\textit{arXiv:2505.04563 [math.CO]}, 2025.

%\bibitem{AK1994}
%K. Aomoto and M. Kita,
%\emph{Theory of hypergeometric functions},
%Springer-Tokyo, 1994 (in Japanese),
%With an appendix by Toshitake Kohno,
%Eng. transl. by K. Iohara, Springer-Verlag, Tokyo,
%2011.

\end{thebibliography}
\end{document}